\documentclass[a4paper,11pt]{article}

\usepackage{amsmath}
\usepackage{amsthm}
\usepackage{amssymb}
\usepackage{cancel}
\usepackage{color}
\usepackage{calrsfs}
\usepackage{varioref}
\usepackage{bbm}
\definecolor{dblue}{rgb}{0.09,0.32,0.44} 
\usepackage[pdfborder={0 0 0}, colorlinks=true, pdfborderstyle={}, linkcolor=dblue, citecolor=dblue, urlcolor=blue]{hyperref}
\usepackage{epsfig}
\usepackage{psfrag}
\usepackage{subfigure}
\usepackage{graphicx}
\usepackage[mathscr,mathcal]{eucal}
\usepackage[shortlabels]{enumitem}
\usepackage{t1enc}
\usepackage[latin2]{inputenc}

\newtheorem {theorem}{Theorem}

\newtheorem {proposition}{Proposition}

\newtheorem* {theorem*}{Theorem}
\newtheorem* {lemma*}{Lemma}
\newtheorem* {corollary*}{Corollary}
\newtheorem* {proposition*}{Proposition}
\newtheorem* {definition*}{Definition}
\newtheorem* {conjecture*}{Conjecture}
\newtheorem* {question*}{Question}

\newtheorem* {theoremkv*} {Theorem KV}
\newtheorem* {corollarykv*} {Corollary KV}
\newtheorem* {theoremrsc1*} {Theorem RSC1}
\newtheorem* {theoremrsc2*} {Theorem RSC2}

\theoremstyle{remark}
\newtheorem* {remark*}{Remark}

\def \N {\mathbb N}

\def \R {\mathbb R}

\def \Z {\mathbb Z}

\def\cB{\mathcal{B}}

\def\cE{\mathcal{E}}
\def\cF{\mathcal{F}}

\def\cH{\mathcal{H}}

\def\cL{\mathcal{L}}

\def\cX{\mathcal{X}}

\def\vareps{\varepsilon}

\newcommand{\probab}[1]{\ensuremath{\mathbf{P}\left(#1\right)}}
\newcommand{\expect}[1]{\ensuremath{\mathbf{E}\left(#1\right)}}

\newcommand{\condexpect}[2]{\ensuremath{\mathbf{E}\left(#1\bigm|#2\right)}}

\newcommand{\probabom}[1]{\ensuremath{\mathbf{P}_{\omega}\left(#1\right)}}
\newcommand{\expectom}[1]{\ensuremath{\mathbf{E}_{\omega}\left(#1\right)}}

\newcommand{\condprobabom}[2]{\ensuremath{\mathbf{P}_{\omega}\left(#1\bigm|#2\right)}}
\newcommand{\condexpectom}[2]{\ensuremath{\mathbf{E}_{\omega}\left(#1\bigm|#2\right)}}

\DeclareMathOperator{\grad}{grad}
\let\div\undefined
\DeclareMathOperator{\div}{div}

\def\Dom{\mathrm{Dom}}

\renewcommand{\d}{\mathrm d}

\newcommand{\abs}[1]{\ensuremath\left|{#1}\right|}
\newcommand{\norm}[1]{\ensuremath\left\|{#1}\right\|}

\def\wh{\widehat}

\newcommand{\wick}[1]{\ensuremath{\,:\! #1 \!:\,}}

\begin{document}

\title{Diffusive and Super-Diffusive Limits for Random Walks and Diffusions with Long Memory
\footnote{ICM-2018 Probability Section talk.
In: {\sl Proceedings of the ICM - 2018 Rio de Janeiro}, Vol. 3 pp 3025-3044, World Scientific 2018
}
}
\author{
{\sc B\'alint T\'oth
\footnote{Supported by EPSRC (UK) Established Career Fellowship EP/P003656/1 and by OTKA (HU) K-109684.}
}
\\
{University of Bristol and R\'enyi Institute, Budapest}
}

\maketitle

\begin{abstract}
\noindent
We survey recent results of normal and anomalous diffusion of two types of random motions with long memory in $\R^d$ or $\Z^d$. The first class consists of random walks on $\Z^d$ in divergence-free random drift field, modelling the motion of a particle suspended in time-stationary incompressible turbulent flow. The second class consists of self-repelling random diffusions, where the diffusing particle is pushed by the negative gradient of its own occupation time measure towards regions less visited in the past. We establish normal diffusion (with square-root-of-time scaling and Gaussian limiting distribution) in three and more dimensions and typically anomalously fast diffusion in low dimensions (typically, one and two). Results are quoted from various papers published between 2012-2018, with some hints to the main ideas of the proofs. No technical details are presented here. 



\end{abstract}

\section{Random walks in divergence-free random drift field}
\label{s: Random walks in divergence-free random drift field}

\subsection{Set-up and notation}
\label{ss: Set-up and notation}

Let $(\Omega, \cF, \pi, \tau_z:z\in\Z^d)$ be a probability space with an ergodic $\Z^d$-action. Denote by  $\cE:=\{k\in\Z^d: |k|=1\}$ the set of possible steps of a nearest-neighbour walk on $\Z^d$, and let $p_k:\Omega\to[0,s^*]$, $k\in\cE$, be bounded measurable functions. These will be the jump rates of the RWRE considered (see \eqref{the walk} below) and assume they are \emph{doubly stochastic}, 
\begin{align}
\label{bistoch}
\sum_{k\in\cE}p_k(\omega)
=
\sum_{k\in\cE}p_{-k}(\tau_k\omega).
\end{align}
Given these,  define the continuous time nearest neighbour random walk $t\mapsto X(t)\in\Z^d$ as a Markov process on $\Z^d$, with  $X(0)=0$ and conditional jump rates
\begin{align}
\label{the walk}
\condprobabom{X(t+dt)= x+k}{X(t)=x} = p_k(\tau_x\omega) dt,
\end{align}
where the subscript $\omega$ denotes that the random walk $X(t)$ is a Markov process on $\Z^d$ \emph{conditionally}, with fixed $\omega\in\Omega$, sampled according to $\pi$. The continuous time setup is for convenience only. Since the jump rates are bounded this is fully equivalent with a discrete time walk. 

We will use the notation $\probabom{\cdot}$ and $\expectom{\cdot}$  for \emph{quenched} probability and expectation. That is: probability and  expectation with respect to the distribution of the random walk $X(t)$, \emph{conditionally, with given fixed environment $\omega$}. The notation $\probab{\cdot}:=\int_\Omega\probabom{\cdot} {\d}\pi(\omega)$ and  $\expect{\cdot}:=\int_\Omega\expectom{\cdot} {\d}\pi(\omega)$ will be used for \emph{annealed} probability and expectation. That is: probability and  expectation with respect to the random walk trajectory $X(t)$ \emph{and} the environment $\omega$, averaged out with  the distribution $\pi$. 

It is well known (and easy to check, see e.g.\ \cite{kozlov-85}) that due to double stochasticity \eqref{bistoch} the annealed set-up is stationary and ergodic in time: the process of the environment as seen from the position of the random walker
\begin{align}
\label{env proc}
\eta(t):=\tau_{X(t)}\omega
\end{align}
is a stationary and ergodic Markov process on $(\Omega, \pi)$ and, consequently, the random walk $t\mapsto X(t)$ will have stationary and ergodic annealed increments.

The local \emph{quenched} drift of the random walk is
\begin{align}
\notag
\condexpectom{dX(t)}{X(t)=x} 
= 
\sum_{k\in\cE} kp_k(\tau_x\omega)dt
=:
\varphi(\tau_x\omega)dt. 
\end{align}
It is convenient to separate the symmetric and skew-symmetric part of the jump rates:  for $k\in\cE$, let $s_k:\Omega\to[0,s^*]$, $v_k:\Omega\to[-s^*,s^*]$,
\begin{align}
\label{s and v}
s_k(\omega):=\frac{p_k(\omega)+p_{-k}(\tau_k\omega)}{2},
&&
v_k(\omega):=\frac{p_k(\omega)-p_{-k}(\tau_k\omega)}{2}.
\end{align}
Note that from the definitions \eqref{s and v} it follows that
\begin{align}
\label{symmetries}
s_k(\omega)-s_{-k}(\tau_k\omega)=0,
&&
v_k(\omega)+v_{-k}(\tau_k\omega)=0.
\end{align}
In addition, the bi-stochasticity condition \eqref{bistoch} is equivalent to 
\begin{align}
\label{divfree}
\sum_{k\in\cE}v_k(\omega)\equiv0, 
\qquad 
\pi\text{-a.s.}
\end{align}
The second identity in \eqref{symmetries} and \eqref{divfree} jointly mean that $\left(v_k(\tau_x\omega)\right)_{k\in\cE, x\in\Z^d}$ is a stationary \emph{sourceless flow} (or, a \emph{di\-ver\-gence-free lattice vector field}) on $\Z^d$. The physical interpretation of the divergence-free condition \eqref{divfree} is that the walk \eqref{the walk} models the motion of a particle suspended in stationary, \emph{incompressible flow}, with thermal noise. 

In order that the walk $t\mapsto X(t)$ have \emph{zero annealed mean drift} we assume that for all $k\in\cE$
\begin{align}
\label{nodrift}
\int_\Omega v_k(\omega) \,{\d}\pi(\omega)
=0.
\end{align}

Our next assumption is an \emph{ellipticity} condition for the symmetric part of the jump rates: there exists another positive constant  $s_*\in(0,s^*]$ such that for $\pi$-almost all $\omega\in\Omega$ and all $k\in\cE$
\begin{align}
\label{ellipt}
s_k(\omega)\ge s_*, 
\quad 
\pi\text{-a.s.}
\end{align}
Note that the ellipticity condition is imposed only on the symmetric part $s_k$ of the jump rates and not on the jump rates $p_k$. It may happen that $\pi( \{ \omega: \min_{k\in\cE}p_k(\omega)=0 \} )>0$, as it is the case in some of the examples given in section \ref{ss: Examples}. 

Finally, we formulate the notorious \emph{$\cH_{-1}$-condition} which plays a key role in diffusive scaling limits. Denote for $i,j=1,\dots,d$, $x\in\Z^d$, $p\in[-\pi,\pi)^d$,  
\begin{align}
\label{covs}
&
C_{ij}(x)
:=
\int_\Omega \varphi_i(\omega)\varphi_j(\tau_x\omega)d\pi(\omega),
&
\wh C_{ij}(p)
:=
\sum_{x\in\Z^d} e^{\sqrt{-1}x\cdot p} C_{ij}(x).
\end{align}
That is: $C_{ij} (x)$ is the covariance matrix of the drift field, and $\wh C_{ij}(p)$ is its Fourier-transform. 

By Bochner's theorem, the Fourier transform $\wh C$ is  positive definite $d\times d$-matrix-valued-measure on $[-\pi,\pi)^d$. The no-drift condition \eqref{nodrift} is equivalent to $\wh C_{ij}(\{0\})=0$, for all $i,j=1,\dots,d$. With slight abuse of notation we denote this measure formally as $\wh C_{ij}(p)dp$ even though it could be not absolutely continuous with respect to Lebesgue. 

The \emph{$\cH_{-1}$-condition} is the following: 
\begin{align}
\label{hcond}
\int_{[-\pi,\pi)^d}
\left(\sum_{j=1}^d(1-\cos p_j) \right)^{-1} \sum_{i=1}^d \wh C_{ii}(p) \, {\d}p <\infty.
\end{align}
This is an \emph{infrared bound} on the correlations of the drift field, $x\mapsto\varphi(\tau_x\omega)\in\R^d$. It implies diffusive upper bound on the annealed variance of the walk and turns out to be a natural sufficient condition for the diffusive scaling limit (that is, CLT for the annealed walk). We'll see further below some other equivalent formulations of the $\cH_{-1}$-condition \eqref{hcond}. Note that the $\cH_{-1}$-condition \eqref{hcond} formally implies the no-drift condition \eqref{nodrift}.

For later reference we state here the closely analogous problem of diffusion in divergence-free random drift field. Let $(\Omega, \cF, \pi, \tau_z: z\in\R^d)$ be now a probability space with an ergodic $\R^d$-action, and $F:\Omega\times \R^d \to \R^d$ a stationary vector field which is $\pi$-almost-surely $C^1$ and divergence-free: 
\begin{align}
\label{divfree 2}
\div F \equiv 0, 
\qquad
\pi\text{-a.s.}
\end{align} 
The diffusion considered is
\begin{align}
\label{sde}
dX(t)=dB(t) + F(X(t))dt. 
\end{align}
The SDE \eqref{sde} has unique strong solution, $\pi$-almost surely. The main question is the same as in the case of the random walk \eqref{the walk}: What is the asymptotic scaling behaviour and scaling limit of $X(t)$, as $t\to\infty$? Under what conditions does the central limit theorem with diffusive scaling and Gaussian limit distribution hold? Although the physical phenomena described by \eqref{the walk}-\eqref{bistoch} and \eqref{sde}-\eqref{divfree 2} are very similar, the technical details of various proofs are not always the same. In particular, PDE methods and techniques used for the diffusion problem \eqref{sde}-\eqref{divfree 2} are not always easily implementable for the lattice problem \eqref{the walk}-\eqref{bistoch}. On the other hand, often restrictive local regularity conditions must be imposed on the diffusion problem \eqref{sde}-\eqref{divfree 2}.

The results reported in this section refer mainly to the random walk problem \eqref{the walk}-\eqref{bistoch}. The diffusion problem \eqref{sde}-\eqref{divfree 2} will be tangentially mentioned in a example in section \ref{sss: Superdiffusive bounds} and in the historical notes of section \ref{ss: Historical notes}. 

\subsection{The infinitesimal generator of the environment process}
\label{ss: The infinitesimal generator of the environment process}

All forthcoming analysis will be done in the Hilbert space 
$\cH
:=
\{ f\in\cL^2(\Omega, \pi): \ 
\int_{\Omega} f(\omega)d\pi(\omega)=0\}.
$
The $\cL^2(\Omega,\pi)$-gradients and Laplacian are bounded operators on $\cH$: 
\begin{align}
\notag
\nabla_k f(\omega)
:=
f(\tau_k\omega)-f(\omega)
&&
\Delta
:=
2\sum_{k\in\cE}\nabla_k
=
-
\sum_{k\in\cE}
\nabla_{-k}\nabla_k.
\end{align}
Note that $\Delta$ is self-adjoint and negative. Thus, the operators $\abs{\Delta}^{1/2}$ and  $\abs{\Delta}^{-1/2}$ are defined in terms of the spectral theorem. The domain of the unbounded operator $\abs{\Delta}^{-1/2}$ is
\begin{align}
\notag
\cH_{-1}:=
\{\phi\in\cH: 
\lim_{\lambda\searrow0}
(\phi, (\lambda I -\Delta)^{-1}\phi)_{\cH}<\infty\}.
\end{align}
The $\cH_{-1}$-condition gets its name from the fact that \eqref{hcond} is equivalent to requesting that for $k\in\cE$,  
\begin{align}
\label{hcond2}
v_k\in\cH_{-1}.
\end{align}

We will also use the multiplication operators $M_k, N_k: \cL^2(\Omega, \pi)\to\cL^2(\Omega, \pi)$, $k\in\cE$, 
\begin{align}
\notag
&
M_kf(\omega):=v_k(\omega)f(\omega),
&&
N_kf(\omega):=\left(s_k(\omega)-s_*\right)f(\omega). 
\end{align}
The following commutation relations are direct consequences of (in fact, equivalent with) \eqref{symmetries} and \eqref{divfree}
\begin{align}
\label{commute}
\sum_{k\in\cE}M_k\nabla_k
=
- 
\sum_{k\in\cE}\nabla_{-k}M_k, 
&&
\sum_{k\in\cE}N_k\nabla_k 
=
\sum_{k\in\cE}\nabla_{-k}N_k
\end{align}
We will denote
\begin{align}
\notag
S
:=
-\frac{s_*}{2}
\Delta
+
\sum_{k\in\cE}N_k\nabla_{k},
=S^*
&&
A
:=
\sum_{k\in\cE}M_k\nabla_k,
=-A^*.
\end{align}
The infinitesimal generator $L$ of the Markovian semigroup $P_t:\cL^2(\Omega, \pi)\to\cL^2(\Omega, \pi)$ of the environment process \eqref{env proc} is 
\begin{align}
\notag
L
=
-
S+A
.
\end{align}
Note that due to ellipticity \eqref{ellipt} and boundedness of the jump rates the (absolute value of the) Laplacian minorizes and majorizes the self-adjoint part of the infinitesimal generator: 
$s_* \abs{\Delta} \le 2 S \le s^* \abs{\Delta}$. The inequalities are meant in operator sense.

\subsection{Helmholtz's theorem and the stream tensor}
\label{ss: Helmholtz's theorem and the stream tensor}

In its most classical form Helmholtz's theorem states that in $\R^3$ (under suitable conditions of moderate increase at infinity) a divergence-free vector field can be realised as the curl (or rotation) of another vector field, called the vector potential. Helmholtz's theorem in our context is the following: 

\begin{proposition}
\label{prop:helmholtz}
Let $v: \Omega\to\R^{\cE}$ be such that $v_k\in\cH$, and assume that \eqref{symmetries} and \eqref{divfree} hold. 

\smallskip
\noindent
(i)
There exists a zero mean, square integrable, antisymmetric tensor cocycle
${H}:\Omega\times\Z^d\to\R^{\cE\times\cE}$, $H_{k,l}(\cdot,x)\in\cH$: 
\begin{align}
\label{Hcocycle}
&
H_{k,l}(\omega,y)-H_{k,l}(\omega,x)
=
H_{k,l}(\tau_x\omega, y-x)-H_{k,l}(\tau_x\omega, 0),
\\
\label{Htensor}
&
H_{l,k}(\omega,x)
=
H_{-k,l}(\omega,x+k)
=
H_{k,-l}(\omega,x+l)
=
-H_{k,l}(\omega,x),
\end{align}
such that
\begin{align}
\label{vcurl}
v_k(\tau_x\omega)=\sum_{l\in\cE} {H}_{k,l}(\omega,x).
\end{align}
The realization of the tensor field ${H}$ is unique up to an additive term $H^0_{k,l}(\omega)$, not depending on $x\in\Z^d$ (but obeying the symmetries \eqref{Htensor}).

\smallskip
\noindent
(ii)
The $\cH_{-1}$-condition \eqref{hcond}/\eqref{hcond2} holds if and only if the cocycle $H$ in (i) is stationary. That is, there exists $h:\Omega\to \R^{\cE\times\cE}$,  with $h_{k,l}\in\cH$, such that
\begin{align}
\label{h is a tensor}
h_{l,k}(\omega)
=
h_{-k,l}(\tau_k\omega)
=
h_{k,-l}(\tau_l\omega)
=
-h_{k,l}(\omega),
\end{align}
and
\begin{align}
\label{v is curl of h}
v_k(\omega)=\sum_{l\in\cE} h_{k,l}(\omega).
\end{align}
The tensor field $H$ is realized as the stationary  lifting of $h$:
$H_{k,l}(\omega,x)=h_{k,l}(\tau_x\omega)$.
\end{proposition}

The fact that $v$ is expressed in \eqref{vcurl} as the curl of the tensor field $H$ having the symmetries \eqref{Htensor}, is essentially the lattice-version of Helmholtz's theorem. Note that \eqref{Htensor} means that the stream tensor field $x\mapsto H(\omega,x)$ is in fact a function of the \emph{oriented plaquettes} of $\Z^d$. In particular, in two-dimensions $x\mapsto H(\omega,x)$ defines a \emph{height function} with stationary increments, on the dual lattice $\Z^{2}+ (1/2,1/2)$, in three-dimensions $x\mapsto H(\omega,x)$ defines an \emph{oriented flow} (that is: a lattice vector field) with stationary increments on the dual lattice $\Z^{3}+ (1/2,1/2,1/2)$. In Helmholtz's theorem, if $d>2$, there is much freedom in the choice of the \emph{gauge} of $H$. The cocycle condition \eqref{Hcocycle} is met by the \emph{Coulomb gauge}, which makes the construction essentially uniquely determined.

In \cite{kozma-toth-17} it was shown that for a RWRE \eqref{the walk} whose environment satisfies conditions \eqref{bistoch}, \eqref{ellipt} and \eqref{hcond} the central limit theorem holds, under diffusive scaling and Gaussian limit with finite and nondegenerate asymptotic covariance, \emph{in probability with respect to the environment}. See Theorem \ref{thm:clt in probability} below. 

In order to obtain the quenched version, that is central limit theorem for the displacement $X(t)$ at late times, with frozen environment, $\pi$-almost surely, we impose a slightly stronger integrability condition on the stream-tensor-field,
\begin{align}
\label{two-plus-eps}
h\in\cL^{2+\vareps} (\Omega, \pi),
\end{align}
for some $\vareps>0$, rather than being merely square integrable. This stronger integrability condition is needed in the proof of quenched tightness of the diffusively scaled displacement $t^{-1/2}X(t)$. We will refer to the $\cH_{-1}$-condition complemented with the stronger integrability assumption \eqref{two-plus-eps} as the \emph{turbo-$\cH_{-1}$-condition}. In \cite{toth-17} the quenched version of the central limit theorem for the displacement of the random walker was proved under the conditions \eqref{bistoch}, \eqref{ellipt} and \eqref{hcond} and \eqref{two-plus-eps}. See Theorem \ref{thm: quenched clt} below.

\subsection{Examples}
\label{ss: Examples}

\noindent
{\bf Bounded stream tensor:}
Let $\left((\chi_{ij}(x))_{1\le i<j\le d}\right)_{x\in\Z^d}$, be a stationary and ergodic (with respect to $x\in\Z^d$) sequence of bounded random variables (say, $\abs{\chi_{ij}(x)}\le1$), and extend them to $i,j\in\{1,\dots,d\}$ as $\chi_{ji}=-\chi_{ij}$, $\chi_{ii}=0$. Define for $k,l\in\cE$, $x\in\Z^d$,
\begin{align}
\notag 
H_{k,l}(x)
:=
(k\cdot e_i) (l\cdot e_j)
\chi_{ij}
( x+  (k\cdot e_i-1)e_i/2  + (l\cdot e_j-1)e_j/2).
\end{align}
(This formula extends the random variables $\chi$ to a tensor field, consistent with the symmetries \eqref{Htensor}). Define $v_k(\omega)$ as in \eqref{vcurl} and let $s_k(\omega)\equiv s_*\ge 2(d-1)$. This is the most general construction of the case when $h\in\cL^\infty(\Omega,\pi)$. In particular, it covers  those cases when the random environment of jump probabilities admits a \emph{bounded cycle representation}, cf. \cite{kozlov-85}, \cite{komorowski-olla-03a}, \cite{deuschel-kosters-08}, \cite{komorowski-landim-olla-12} (chapter 3.3). Due to Proposition \ref{prop:helmholtz} the $\cH_{-1}$-condition \eqref{hcond}/\eqref{hcond2} holds. 

\smallskip
\noindent
{\bf Randomly oriented Manhattan lattice:}
Let $u_i(y)$, $i\in\{1,\dots,d\}$, $y\in \Z^{d-1}$, be i.i.d. random variables with the common distribution $\probab{u=\pm1}=1/2$, and define for $x\in\Z^d$ and $k\in \cE$
\begin{align}
\notag
v_k(\tau_x\omega)
:=
\sum_{i=1}^d (k\cdot e_i) u_i(x_1,\dots, x_{i-1}, \cancel{x_i}, x_{i+1}, \dots, x_d),
\end{align}
One can easily compute the covariances \eqref{covs}:
$C_{ij}(x)=\delta_{i,j} \prod_{i'\not=i}\delta_{x_{i'},0}$, and their Fourier transforms
$\wh C_{ij}(p) = \delta_{i,j} \delta(p_i)$.
From here it follows that in this particular model the $\cH_{-1}$-condition fails robustly (with power law divergence in \eqref{hcond}) in $d=2$, fails marginally (with logarithmic divergence in \eqref{hcond}) in $d=3$,  and holds if $d\ge 4$. 

\smallskip
\noindent
{\bf $\binom{2d}{d}$-vertex models on $\Z^d$:}
Let $\Omega$ be the set of all possible orientations of the edges of $\Z^d$ with the constraint that at all vertices the number of edges oriented towards the site is equal to the number of edges oriented away, $d$ out of $2d$. In this way, locally at every vertex $\binom{2d}{d}$ configurations of orientations are possible and there is a very rigid constraint on the configurations. $\Omega$ is a compact metric space and the group of translations $\tau_z:\Omega\to\Omega$, $z\in\Z^d$, acts naturally on it. Let, for $k\in\cE$,  $v_k(\omega)=\pm1$ be the orientation of the edge $\overline{0k}$ in the configuration $\omega\in\Omega$, and $p_k(\omega)=1+v_k(\omega)$. Any translation invariant ergodic measure $\pi$ on $\Omega$ realizes a model of our RWRE. The most natural choice is the one when $\pi$ is the unique weak limit of the uniform distribution of the allowed $\binom{2d}{d}$-vertex configurations on the discrete torus $(-L, L]\times\dots\times (-L,L]$, with periodic boundary conditions, as $L\to\infty$. In 2-dimensions this is the notorious (uniform) six-vertex model. In this case - in 2-dimensions - the $\cH_{-1}$-condition fails: the integral in \eqref{hcond} is logarithmically divergent. 

\subsection{Scaling limits}
\label{ss: Sscaling limits}

\subsubsection{Central limit theorem in probability w.r.t. the environment under the $\cH_{-1}$-condition}
\label{sss: Annealed central limit theorem}

\begin{theorem}
\label{thm:clt in probability}
{\rm (Source: \cite{kozma-toth-17})}
Conditions \eqref{bistoch}, \eqref{ellipt}, \eqref{hcond} are assumed. 
The asymptotic annealed covariance matrix
\begin{align}
\label{limit covariance}
(\sigma^2)_{ij}
:=
\lim_{t\to\infty} t^{-1}\expect{X_i(t)X_j(t)}
\end{align}
exists, and it is finite and non-degenerate.
For any bounded and continuous function $f:\R^d\to\R$, 
\begin{align}
\label{clt in prob}
\lim_{T\to\infty}
\int_\Omega
\abs{
\expectom{f(T^{-1/2}X(T))}
-
\int_{\R^d} \frac{ e^{-\frac {\abs{y}^2} {2} }}{(2\pi)^{\frac{d}{2}}} f(\sigma^{-1}y) dy
}
d\pi(\omega)=0.
\end{align}
\end{theorem}

Theorem \ref{thm:clt in probability} is proved in \cite{kozma-toth-17}, and, weak convergence in the sense of \eqref{clt in prob} of all finite dimensional marginal distributions of $t\mapsto T^{-\frac12} X(Tt)$, as $T\to\infty$, to those of a $d$-dimensional Brownian motion with covariance $\sigma^2$ is established.  We sketch the main points. Start with a most natural martingale decomposition 
\begin{align}
\label{martingale decomposition}
X(t)
=
\left\{
X(t)-
\int_{0}^{t} \varphi(\eta(s))ds
\right\}
+
\int_{0}^{t} \varphi(\eta(s))ds
=:
M(t)+I(t). 
\end{align}
In this decomposition $M(t)$ is clearly a square integrable martingale with stationary and ergodic annealed increments. The main issue is an efficient martingale approximation of the term $I(t)$, \`a la Kipnis-Varadhan. 

We rely on the following general result on Kipnis-Varadhan type of martingale approximation. Let $\eta(t)$ be a stationary and ergodic Markov process on the probability space $(\Omega,\pi)$, and $L$ be the infinitesimal generator of its Markovian semigroup acting on $\cL^2(\Omega,\pi)$. Denote $S:=-(L+L^*)/2$, $A:=(L-L^*)/2$ and assume that the symmetric part $S$ is minorised and majorized by a postitive operator $D\ge0$: $s_* D \le S\le s^* D$, with $0<s_*\le s^*<\infty$.  
Further, denote
\begin{align}
\notag
B_\lambda:=(\lambda I + D)^{-1/2} A (\lambda I + D)^{-1/2}.
\end{align}

\begin{theorem}
\label{thm:rsc}
{\rm (Source: \cite{horvath-toth-veto-12a}, \cite{kozma-toth-17})}
Assume that there exist a dense subspace $\cB\subseteq \cL^2(\Omega,\pi)$ and a linear operator $B:\cB\to \cL^2(\Omega,\pi)$ which is essentially skew-self-adjoint on the core $\cB$  and such that for any $\varphi\in \cB$ there exists a sequence $\varphi_\lambda\in\cL^2(\Omega,\pi)$ such that
\begin{align}
\notag
\lim_{\lambda\to 0}\norm{\varphi_\lambda-\varphi}=0.
\qquad
\text{and}
\qquad
\lim_{\lambda\to 0}\norm{B_\lambda\varphi_\lambda-B\varphi}=0.
\end{align}
Then for any $f\in\Dom(D^{-1/2})$ there exists a martingale $M_f(t)$ (on the probability space and with respect to the filtration of the Markov process $t\mapsto\eta_t$) such that 
\begin{align}
\notag
\lim_{t\to\infty} t^{-1}\int_\Omega \expectom{\abs{\int_0^t f(\eta(s)) ds - M_f(t)}^2} =0.
\end{align}
\end{theorem}

In plain and informal words: if the operator $B=D^{-1/2} A D^{-1/2}$ makes sense as a densely defined  unbounded \emph{skew-self-adjoint} operator then integrals along the Markov process trajectory of functions in $\cH_{-1}\subset\cL^2(\Omega, \pi)$ defined in terms of the positive operator $D$ are efficiently approximated with martingales, \`a la Kipnis-Varadhan. As shown in \cite{horvath-toth-veto-12a}, the condition of Theorem \ref{thm:rsc} is weaker than the graded sector condition of \cite{sethuraman-varadhan-yau-00} (which in turn is weaker than the strong sector condition of \cite{varadhan-95}).

In our particular case define $B:\cH_{-1}\to\cH$ as
\begin{align}
\label{Bop def}
B
:=
-
\sum_{l\in\cE}\abs{\Delta}^{-1/2}\nabla_{-l} M_l \abs{\Delta}^{-1/2}.
\end{align}
(Note that the operators $\abs{\Delta}^{-1/2}\nabla_{-l}$, $l\in\cE$ are bounded.)
From the commutation relations \eqref{commute} it follows that the  operator $B$ is \emph{skew-symmetric} on the dense subspace $\cB:=\cH_{-1}$. 
It is not difficult to show that if $h\in\cL^\infty$ then $B$ is a bounded operator and thus, its skew-self-adjointness drops out for free. (This is essentially the same as Varadhan's strong sector condition, cf.\cite{varadhan-95}) On the other hand, if $h\notin\cL^\infty$ then $B$ is genuinely unbounded and proving its (essential) skew-self-adjointness is far from trivial. 

The main technical result in \cite{kozma-toth-17} is the proof of the fact that $B$ is in fact \emph{essentially skew-self-adjoint} on $\cH_{-1}$. By applying von Neumann's criterion for (skew-)self-adjointness this boils down to proving that the lattice PDE
\begin{align}
\label{lattice pde}
\Delta\Psi(\cdot, \omega) +  V(\cdot, \omega) \cdot \nabla \Psi(\cdot, \omega) = 0,
\end{align}
does not have a non-trivial cocycle solution $\Psi(x,\omega)$. Here now $\nabla$ and $\Delta$ denote the lattice gradient, respectively, the lattice Laplacian, $V(x,\omega)=v(\tau_x\omega)$ and $\Psi(x, \omega)$ is a zero mean cocycle to be determined. 

\subsubsection{Quenched central limit theorem under the turbo-$\cH_{-1}$-condition}
\label{sss: Quenched central limit theorem}

\begin{theorem}
\label{thm: quenched clt}
{\rm (Source: \cite{toth-17})}
Conditions \eqref{bistoch}, \eqref{ellipt}, \eqref{hcond}, and \eqref{two-plus-eps} are assumed. 
For any bounded continuous function $f:\R^d\to\R$, 
\begin{align}
\notag
\lim_{T\to\infty}
\expectom{f(T^{-1/2}X(T))}
=
\int_{\R^d} \frac{ e^{-\frac {\abs{y}^2} {2} }}{(2\pi)^{\frac{d}{2}}} f(\sigma^{-1}y) dy
,
&&
\pi\text{-a.s.}
\end{align}
with the non-degenerate covariance matrix $\sigma^2$ given in \eqref{limit covariance}.
\end{theorem}

Theorem \ref{thm: quenched clt} is proved in \cite{toth-17}, and as in the case of Theorem \ref{thm:clt in probability}, the weak convergence of all finite dimensional distributions follows. The proof consists of three major steps: 
(1) 
Proof of \emph{quenched tightness} of the scaled displacement $t^{-1/2} X(t)$, as $t\to\infty$.
(2)
Construction of the \emph{harmonic coordinates} for the walk. 
(3)
Efficient estimate of the discrepancy between the actual position of the walker and the approximating harmonic coordinates. 

\smallskip
\noindent
{\bf Quenched tightness of $t^{-1/2} X(t)$:}
\ \ \ 
The proof relies on adapting Nash's moment bound on reversible diffusions with strictly elliptic and bounded dispersion coefficients, cf. \cite{nash-58}, to this type of non-reversible setup. The extension in the case of $h\in\cL^\infty(\Omega, \pi)$ is essentially straightforward, following \cite{osada-83}. (Though, adaptation to the lattice walk case needs some attention.) The extension to $h\in \cL^{2+\vareps}(\Omega, \pi)$ is tricky. An integration over time and the Chacon-Ornstein ergodic theorem help. Full details can be found in \cite{toth-17}.  This is the only part of the proof where the stronger integrability condition \eqref{two-plus-eps} is used. 

\smallskip
\noindent
{\bf Harmonic coordinates:}
\ \ \ 
The idea of harmonic coordinates for random walks in random environments originates in the classical works 
\cite{kozlov-80}, \cite{papanicolaou-varadhan-81}, \cite{osada-83}, \cite{kozlov-85}. Since then it had been widely used in proving quenched central limit theorems, mostly for random walks among random conductances. That is: in time reversible cases. See, however, \cite{deuschel-kosters-08} for a non-reversible application. The idea is very natural: find an $\R^d$-valued $\cL^2(\Omega,\pi)$, zero mean  random cocycle $\Theta(x, \omega)$, such that 
\begin{align}
\label{harmcoordeq}
\sum_{k\in\cE}p_k(\tau_x\omega)\left(k+\Theta(\omega,x+k) - \Theta(\omega,x)\right)=0, 
\quad
\pi\text{-a.s.}
\end{align}
If there exists a solution $\Theta$ to the equation \eqref{harmcoordeq} then the process 
$t\mapsto Y(t):=X(t) + \Theta(\omega, X(t))$
is a quenched martingale (that is, a martingale in its own filtration, with the environment $\omega\in\Omega$ frozen). It turns out that equation \eqref{harmcoordeq} is equivalent with the following equation in $\cH$:
\begin{align}
\label{eqforchi}
(I+ B^*)\chi = \abs{\Delta}^{-1/2}\varphi, 
\end{align}
where $\varphi\in\cH_{-1}$ is given and $\chi$ is to be determined. The operator $B^*$ on the left hand side is exactly the adjoint of $B$ from \eqref{Bop def}. Since it was proved that the operator $B$ is skew-self-adjoint it follows that $I+B^*$ is invertible and thus equation \eqref{eqforchi} has a unique solution. As a consequence, equation \eqref{harmcoordeq} also has a unique solution $\Theta$ which is an $\R^d$-valued cocycle, as required.

Once the harmonic coordinates are constructed the quenched central limit theorem for $t^{-1/2}Y(t)$ drops out via the martingale CLT, using ergodicity of the environment process \eqref{env proc}. 

\smallskip
\noindent
{\bf Error bound:}
\ \ \ 
In order to obtain the quenched CLT for the scaled displacement $t^{-1/2}X(t)$ it remains to prove that for all $\delta>0$ and  $\pi$-almost all $\omega\in\Omega$, 
\begin{align}
\notag
\lim_{t\to\infty}\probabom{\abs{\Theta( X(t))}>\delta \sqrt{t}}=0.
\end{align}
The key ingredients of this are the a priori quenched tightness of $t^{-1/2}X(t)$ proved in the first main step, and a soft but nevertheless useful ergodic theorem for cocycles: Let $\Omega\times\Z^d\ni x\mapsto \Psi(\omega,x)\in\R$ be a zero-mean $\cL^2$-cocycle. Then
\begin{align}
\notag
\lim_{N\to\infty} N^{-(d+1)} \sum_{\abs{x}\le N} \abs{\Psi(x)}
=0, 
\qquad
\pi\text{-a.s.}
\end{align}
In 1-dimension this is a direct consequence of Birkhoff's ergodic theorem. For $d>1$, however, the multidimensional \emph{unconditional ergodic theorem} is invoked. See \cite{toth-17} for the proof.

\subsubsection{Superdiffusive bounds when the $\cH_{-1}$-condition fails}
\label{sss: Superdiffusive bounds}

If the $\cH_{-1}$-condition \eqref{hcond}/\eqref{hcond2} fails, or, equivalently, the conditions of part ($ii$) of Proposition \ref{prop:helmholtz} don't hold, then there is no a priori upper bound on $t^{-1}\expect{\abs{X(t)}^2}$ and superdiffusive behaviour is expected. There is no general statement like this, but there are particular interesting cases studied. The fully worked out cases are, however, continuous-space diffusions on $\R^d$ rather than random walks on $\Z^d$. Since $d=2$ is the most interesting we only mention the following two-dimensional example.

Let $x\mapsto F(x)$ be a stationary Gaussian random vector field with covariances $K_{ij}(x):=\expect{F_i(0)F_j(x)}$ as follows
\begin{align}
\label{covarF}
K_{ij}(x)
=
\left(
\partial^2_{ij}
-
\delta_{i,j}(\partial^2_{11}+\partial^2_{22})
\right) 
V * G (x), 
&&
\wh K_{ij}(p)
&=
\left(
\delta_{i,j}
-
\frac{p_i p_j}{\abs{p}^2} 
\right)
\wh V(p),
\end{align}
where $G(x)=\log\abs{x}$ is the two-dimensional (Laplacian) Green function and $V:\R^2\to\R_+$ is a $C^\infty$ approximate identity with fast decay and positive Fourier transform, $\wh V(p)>0$. A good concrete choice could be $V(x)=(2\pi \sigma^2)^{-{d}/{2}}\exp\{-{\abs{x}^2}/{(2\sigma^2)}\}$, with some $\sigma\in(0,\infty)$. In plain words: $F$ is the rotation (curl) of the two-dimensional Gaussian free field, locally mollified by convolving with the convolution-square-root of $V$. As a rotation, the vector field $F$ is divergence-free, cf. \eqref{divfree 2}. Define the diffusion in this random drift field: $t\mapsto X(t)\in\R^2$ as the unique strong solution of the SDE \eqref{sde}. 

From \eqref{covarF} it appears that the $\cH_{-1}$-condition fails marginally: the integral on the right hand side of \eqref{hcond} diverges logarithmically. In \cite{toth-valko-12} superdiffusive bounds are proved for this diffusion in the rotation field of the two-dimensional Gaussian free field, which look formally very similar to  \eqref{d=2 superdiff} in section \ref{ss: Scaling limits 2} below. This extends earlier results (with power-law divergences)  of \cite{komorowski-olla-02} to the marginal case of the two-dimensional Gaussian free field (with logarithmic divergences).  The random walk on the six-vertex model (see the third example in section \ref{ss: Examples}) behaves phenomenologically similarly, but its superdiffusivity is not yet proved. Applying the same methods as in \cite{toth-valko-12} we obtain, however, superdiffusive bounds for the variance of $X(t)$ for the random walk on the randomly oriented Manhattan lattice (second example in section \ref{ss: Examples}) in dimensions $d=2$ (robust, power law) and $d=3$ (marginal, logarithmic), cf. \cite{ledger-toth-valko-17}. 

\subsection{Historical notes}
\label{ss: Historical notes}

The problems of scaling limit of diffusions in divergence-free random drift field \eqref{sde}-\eqref{divfree 2} and that of the random walks in doubly stochastic random environment \eqref{the walk}-\eqref{bistoch} are closely related. Although the physical phenomena modelled are very similar (tracer motion along the drift lines of incompressible turbulent flow), the technical details of various proofs are not always the same. In particular, PDE methods and techniques used for the diffusion problem \eqref{sde}-\eqref{divfree 2} are not always easily implementable for the lattice problem \eqref{the walk}-\eqref{bistoch}. In the following list we give a summary of the main stations in the probability literature along the almost forty years history of the subject. The list is far from complete and contains only the probability results. See also the bibliographical notes of chapters 3 and 11 of \cite{komorowski-landim-olla-12}.

\smallskip
\noindent
1979: 
Kozlov, respectively, Papanicolaou and Varadhan, independently and in parallel formulate the problem of scaling limits of diffusions in stationary random environment and prove the first CLT for the self-adjoint case under strong ellipticity condition, \cite{kozlov-80}, \cite{papanicolaou-varadhan-81}.

\smallskip
\noindent
1983: 
Osada proves quebched CLT for the diffusion \eqref{sde} in divergence-free drift field \eqref{divfree 2}, when the the stream tensor is bounded, \cite{osada-83}.  

\smallskip
\noindent
1985:
Kozlov formulates the problem of random walk in doubly stochasic random environment \eqref{bistoch}-\eqref{the walk}. An annealed CLT is stated for the case when the jump probabilities $((p_k(\tau_x\omega))_{k\in\cE})_{x\in\Z^d}$ are \emph{finitely dependent}, \cite{kozlov-85}. Double stochasticity \eqref{bistoch} and finite dependence of $(p(\tau_x\omega))_{x\in\Z^d}$,  jointly are rather restrictive conditions. Natural examples are provided by a Bernoulli soup of bounded cycles. 

\smallskip
\noindent
1988: 
Oelschl\"ager proves annealed invariance principle for the diffusion problem \eqref{sde}-\eqref{divfree 2}, under the optimal $\cH_{-1}$-condition, and local regularity condition imposed on the drift field. 

\smallskip
\noindent
1996: 
Fannjiang and Papanicolaou consider the homogenisation for the pa\-ra\-bo\-lic problem corresponding to \eqref{sde}-\eqref{divfree 2}, under $\cH_{-1}$-condition, \cite{fannjiang-papanicolaou-96}.

\smallskip
\noindent
1997:  
Fannjiang and Komorowski prove quenched invariance principle  for the diffusion \eqref{sde}-\eqref{divfree 2}, under the condition that $h\in\cL^p$, with $p>d$, \cite{fannjiang-komorowski-97}

\smallskip
\noindent
2003: 
Komorowski and Olla prove annealed CLT for the random walk \eqref{the walk}-\eqref{bistoch} when $h\in\cL^\infty$, by applying Varadhan's strong sector condition, \cite{komorowski-olla-03a}.

\smallskip
\noindent
2008: 
Deuschel and K\"osters prove quenched CLT for the random walk \eqref{the walk}-\eqref{bistoch} when the jump probabilities admit a bounded cyclic representation, \cite{deuschel-kosters-08}. This condition implies $h\in\cL^\infty$ and thus the strong sector condition.

\smallskip
\noindent
2012: 
Komorowski, Landim and Olla publish the proof of the CLT for the random walk problem \eqref{the walk}-\eqref{bistoch} in the case when $h\in\cL^p$, $p>d$, \cite{komorowski-landim-olla-12}. 

\smallskip
\noindent
2017: 
Annealed CLT for the doubly stochastic RWRE problem \eqref{the walk}-\eqref{bistoch} under the optimal $\cH_{-1}$-condition is proved in \cite{kozma-toth-17}. Quenched CLT under the additional integrability condition  \eqref{two-plus-eps} is proved in \cite{toth-17}.

\section{Self-repelling Brownian polymers}
\label{s: Self-repelling Brownian polymers}

\subsection{Set-up and notation}
\label{ss: Set-up and notation 2}

We consider a self-repelling random process $t\mapsto X(t)\in\R^d$ which is pushed by the negative gradient of its own occupation time measure, towards regions less visited in the past. In this order, fix an approximate identity $V:\R^d\to\R_+$, which is infinitely differentiable, decays exponentially fast at infinity, and is of \emph{positive type}:
\begin{align}
\label{positive type}
\wh V(p)
:=
(2\pi)^{-d/2} \int_{\R^d} e^{i p\cdot x} V(x) d x
\ge 0.
\end{align}
As an example, take $V(x)=(2\pi \sigma^2)^{-{d}/{2}}e^{-{\abs{x}^2}/{(2\sigma^2})}$, $\wh V(p) = e^{-\sigma^2\abs{p}^2/2}$.

Let $X(t)$ be the unique strong solution of the stochastic differential equation 
\begin{align}
\label{sde2}
d X(t)= d B(t)- \left(\int_0^t \grad V (X(t)-X(u)) du\right) dt.
\end{align}
Denoting by $\ell(t,\cdot)$ the \emph{occupation time measure} of the process $X$, 
\begin{align}
\label{local time}
\ell(t, A):=\abs{\{0<s\le t: X(s)\in A\}}, 
\end{align}
where $A\subset\R^d$ is any measurable set, the SDE  \eqref{sde2} is written in the alternative form 
\begin{align}
\label{sde3}
d X(t)= d B(t) - \grad (V*\ell(t,\cdot))(X(t))d t,
\end{align}
which is more suggestive regarding the nature of the process $t\mapsto X(t)$: it is indeed driven by the negative gradient of an appropriate local regularization of its occupation time measure (local time). Following the terminology of the related probability literature we will refer to the process $X(t)$ defined in \eqref{sde2}/\eqref{sde3} as the \emph{self-repelling Brownian polymer}. The main question is: What is the long time asymptotic behaviour or $X$? How does the self-repulsion of the trajectory influence the long time scaling? The problem arose essentially in parallel, but unrelated, in the physics (random walk versions) and probability (diffusion versions) literature, cf. \cite{amit-parisi-peliti-83}, \cite{obukhov-peliti-83}, \cite{peliti-pietronero-87}, \cite{norris-rogers-williams-87}, \cite{durrett-rogers-92}, \cite{cranston-mountford-96}. Mathematically non-rigorous, nevertheless strong and compelling scaling arguments appearing in physics papers  \cite{amit-parisi-peliti-83}, \cite{obukhov-peliti-83}, \cite{peliti-pietronero-87} convincingly suggest the following dimension dependent asymptotic scaling behaviour, as $t\to\infty$: 
\begin{align}
\label{scaling conj}
X(t)
\sim
\begin{cases}
t^{2/3} & \text{ in } d=1, \text{ with non-Gaussian scaling limit},
\\
t^{1/2} (\log t)^{\gamma} & \text{ in } d=2, \text{ with Gaussian scaling limit}, 
\\
t^{1/2} & \text{ in } d\ge3, \text{ with Gaussian scaling limit.}
\end{cases}
\end{align}
In $d=2$, the value of the exponent $\gamma\in(0,1/2]$ in the logarithmic correction is disputed. However, there is good scaling reason to expect $\gamma=\frac14$. In the following sections we are going to present the mathematically rigorous results related to the conjectures of \eqref{scaling conj}.

In one-dimension, for some particular nearest-neighbour lattice walk versions of the self-repelling motion, the conjecture (in the first line of) \eqref{scaling conj} is fully settled. In \cite{toth-95} a limit theorem is proved for $t^{-2/3}X(t)$, with an intricate non-Gaussian limit distribution, believed to be universally valid for the 1-d cases. In \cite{toth-werner-98} the presumed scaling limit process $t\mapsto \cX(t)$ is constructed and fully analysed. In this note we are not going to cover those older results. 

\subsection{The environment process}
\label{s: The environment process}

Let 
\begin{align}
\notag
\Omega 
:= 
\left\{
\omega\in C^{\infty}(\R^d\to\R^d)
\,:\, 
\partial_k\omega_l=\partial_l\omega_k,
\ \
\norm{\omega}_{k,m,r}<\infty, 
\right\}, 
\end{align}
where the seminorms are, for $k\in\{1,\dots,d\}$, $m\in\N^d$, $r\in\N$,
\begin{align}
\notag
\norm{\omega}_{k,m,r} 
:= 
\sup_{x\in\R^d} 
(1+\abs{x})^{-1/r} \abs{\partial^{m_1+\dots+m_d}_{m_1,\dots,m_d}\omega_k(x)}.
\end{align}
I plain words: $\Omega$ is the space of gradient vector fields on $\R^d$, with all partial derivatives increasing slower than any power of $\abs{x}$. 

It turns out that the process $t\mapsto \eta(t, \cdot) \in \Omega$, 
\begin{align}
\notag
\eta(t,x)
:=
\grad (V*\ell (t,\cdot))(X(t)+x)
\end{align}
is a Markov process on the state space $\Omega$ with almost surely continuous sample paths. We allow for an initial profile $\eta(0,\cdot)\in\Omega$. (This means an initial signed measure $\ell(t,\cdot)$ in \eqref{local time}.) The finite dimensional non-Markovian process $t\mapsto X(t)\in\R^d$ is traded for the infinite dimensional Markov process $t\mapsto \eta(t)\in\Omega$. 

Next we define a Gaussian probability measure on $\Omega$: the distribution of the gradient of the Gaussian free field on $\R^d$, locally regularised by convolving with the convolution-square-root of $V$. This is the point where positive definiteness \eqref{positive type} of the self-interaction potential $V$ is essential. Let $\pi$ be the Gaussian measure on $\Omega$ with zero mean and covariances $K_{kl} (y-x):=\int_\Omega \omega_k(x) \omega_l (y) d\pi(\omega)$, 
\begin{align}
\label{covarom}
K_{kl} (x), 
=
-\partial_{kl} V* G (x), 
&&
\wh K_{kl}(p)= \frac{p_kp_l}{{\abs{p}}^2} \wh V(p).
\end{align}
where $G:\R^d\setminus\{0\} \mapsto \R$ is the (Laplacian) Green function.

The group of translations $\tau_z:\Omega\to\Omega$, $z\in\R^d$, acts naturally as $\tau_z\omega(x):=\omega(z+x)$, and $(\Omega, \cF, \pi, \tau_z: z\in\R^d)$ is ergodic.   

\begin{proposition}
\label{prop:stationary measure}
{\rm (Source: \cite{tarres-toth-valko-12}, \cite{horvath-toth-veto-12b})}
The Gaussian probability measure $\pi$, with zero mean and covariances \eqref{covarom} is stationary and ergodic for the Markov process $\eta(t)$. 
\end{proposition}

This fact is consequence of the harmony between the two mechanisms driving the process $\eta$: diffusion pushed by $\eta(t,0)$ and building up (gradient of) local time. Proposition \ref{prop:stationary measure} has two different proofs. In \cite{tarres-toth-valko-12} it is proved through careful It\^o-calculus. In \cite{horvath-toth-veto-12b} a functional analytic proof is presented. 

Tilted Gaussian measures with nonzero constant expectation and the same covariances as in \eqref{covarom} are also stationary and ergodic. We are not considering them because they result in ballistic behaviour (that is: nonzero overall speed) of the motion $X(t)$. We think (though, don't prove) that in $d=1,2$ these are the only stationary and ergodic probability measures for the Markov process $\eta(t)$. In $d\ge3$, however, other stationary distributions of totally different character do exist. 

The forthcoming results are all valid in this stationary regime. That is: the initial $\eta(0,\cdot)$ is sampled according to the distribution $\pi$. As in \eqref{martingale decomposition}, the displacement of the random walker $X(t)$ will be decomposed as sum of a martingale with stationary and ergodic increments and its compensator
\begin{align}
\label{martingale decomposition 2}
X(t)
= 
\left\{
X(t)- \int_0^t \varphi(\eta(s)) ds
\right\}
+
\int_0^t \varphi(\eta(s)) ds
=: M(t) + I(t), 
\end{align}
where now $\varphi:\Omega\to\R^d$, $\varphi_l(\omega)=\omega_l(0)$.
The first term, $M(t)$, in \eqref{martingale decomposition 2} is a square integrable martingale with stationary and ergodic increments (on the probability space and with respect to the natural filtration of the Markov process $t\mapsto \eta(t)$). So, that part is well understood from start: it is diffusive and the martingale central limit theorem applies to it. In one and two dimensions \emph{superdiffusive lower bounds} have been proved for the second term, $I(t)$, on the right hand side of \eqref{martingale decomposition 2}, cf. \cite{tarres-toth-valko-12}, \cite{toth-valko-12} and Theorem \ref{thm:srbp superdiff} below. On the other hand, in three and more dimensions an efficient martingale approximation \`a la Kipnis-Varadhan holds for the compensator term, $I(t)$, on the right hand side of \eqref{martingale decomposition 2}, cf. \cite{horvath-toth-veto-12b} and Theorem \ref{thm:srbp clt}.

\subsection{The infinitesimal generator of the environment process}
\label{ss: The infinitesimal generator of the environment process 2}

All computations will be performed in the Hilbert space $\cH:=\{f\in \cL^2(\Omega, \pi): \int_\Omega f d\pi=0\}$. Scalar product in the Hilbert space $\cH$ will be denoted $\langle\cdot,\cdot\rangle$. This is a \emph{Gaussian Hilbert space} with its natural grading: 
\begin{align}
\label{grading}
\cH
=
\overline{\bigoplus_{n=1}^\infty \cH_n},
\end{align}
where $\cH_n$ is the subspace spanned by the $n$-fold Wick products \\ $\wick{\omega_{k_1}(x_1),\dots, \omega_{k_n}(x_n)}$, $k_j\in\{1,\dots,d\}$, $x_j\in\R^d$. 

The shift operators $U_z: \cH\to\cH$, $U_zf(\omega):=f(\tau_z\omega)$, $z\in\R^d$ form a unitary representation of $\R^d$. Denote by $\nabla_k$, $k\in\{1,\dots,d\}$ the infinitesimal generators: $U_z=e^{\sum_{k=1}^d z_k \nabla_k}$ and 
$\Delta:=\sum_{k=1}^d \nabla_k^2 = -\sum_{k=1}^d \nabla^*_k\nabla_k$. 
Note that the shift operators and all operators derived from them (e.g. $\nabla_k$, $\Delta$) preserve the grading \eqref{grading}. 

We will also use the \emph{creation and annihilation operators} $a^*_l: \cH_n\to\cH_{n+1}$, $a_l: \cH_n\to\cH_{n-1}$ defined on Wick monomials as follows and extended by linearity. 
\begin{align}
\notag
a^{*}_l
\wick{ \omega_{k_1}(x_1),..., \omega_{k_n}(x_n)}
&=
\wick{\omega_l(0), \omega_{k_1}(x_1), ..., \omega_{k_n}(x_n)}
\\
\notag
a_l
\wick{\omega_{k_1}(x_1), ..., \omega_{k_n}(x_n)}
&=
\sum_{m=1}^n 
K_{lk_m}(x_m)
\wick{\omega_{k_1}(x_1), ..., \cancel{\omega_{k_m}(x_m)},...,\omega_{k_n}(x_n)}
\end{align}
As suggested by notation the operators $a_l$ and $a_l^*$ are adjoints of each other and restricted to any finite grade they are bounded: 
\begin{align}
\notag
\norm{a_l}_{{\cH_{n}}\to\cH_{n-1}}
=
\left(\int_{R^d} {\abs{p}^{-2}} {p_l^2}\wh V(p) dp\right)^{1/2}
\sqrt{n}.
\end{align}
The infinitesimal generator $L$ of the semigroup of the Markov process $\eta(t)$, acting on $\cH$, $P_tf(\omega):=\condexpect{f(\eta(t))}{\eta(0)=\omega}$ is expressed in terms of the operators introduced above, as
\begin{align}
\label{infgen srbp}
L=-\frac12 \Delta 
+ 
\sum_{l=1}^d a_l^* \nabla_l
+ 
\sum_{l=1}^d \nabla_l a_l
=
-S+A_++A_-.
\end{align}
The proof of this form of the infinitesimal generator relies -- beside usual manipulations (integration by parts, etc.) -- on \emph{directional derivative} identities in the Gaussian Hilbert space $\cH$ (that is:  elements of Malliavin calculus). For details see \cite{tarres-toth-valko-12}, \cite{horvath-toth-veto-12b}. The notation indicates that $A_+:\cH_n\to\cH_{n+1}$ while $A_-:\cH_n\to\cH_{n+1}$, and clearly, $S=S^*\ge0$, $A_+=-A_-^*$.

\subsection{Scaling limits}
\label{ss: Scaling limits 2}

\subsubsection{Superdiffusive bounds in $d=1$ and $d=2$}
\label{sss: Superdiffusive bounds in d=1 and d=2}

Let, for $\lambda>0$
\begin{align}
\notag
\wh E(\lambda)
:=
\int_0^\infty e^{-\lambda t}\expect{\abs{X(t)}^2} dt. 
\end{align}

\begin{theorem}
\label{thm:srbp superdiff}
{\rm (Source: \cite{tarres-toth-valko-12}, \cite{toth-valko-12})}
In $d\le2$, the following bounds hold, with some constants $0<c<C<\infty$, as $\lambda\to 0$ : 
\begin{align}
\notag
&
\text{in }d=1:
&&
\hskip15mm
c \lambda^{-\frac94}
<
\wh E(\lambda)
<
C \lambda^{-\frac52},
&
\\
\notag
&
\text{in }d=2:
&&
c \lambda^{-2} \log\abs{\log \lambda}
<
\wh E(\lambda)
<
C \lambda^{-2} \log\abs{\lambda}. 
&
\end{align}
\end{theorem}

\noindent
With some minimal regularity assumtion on the asymptotic behaviour of $t\mapsto\expect{\abs{X(t)}^2}$, as $t\to\infty$, the bounds in Theorem \ref{thm:srbp superdiff} imply bounds on their \emph{C\'esaro means}, 
\begin{align}
\notag
&
\text{in } d=1:
&&
\hskip11mm
c t^{\frac54} 
<
\frac1t
\int_0^t\expect{\abs{X(s)}^2} ds
<
Ct^{\frac32},
&
\\
\label{d=2 superdiff}
&
\text{in } d=2:
&&
c t \log \log t
<
\frac1t
\int_0^t\expect{\abs{X(s)}^2} ds
<
C t \log t. 
&
\end{align}
With some extra work the upper bounds can be improved to hold without C\'esaro averaging, see e.g. \cite{landim-olla-yau-98}. However, the loweer bounds are the more interesting here. The bounds are consistent with but don't quite match the asymptotic behaviour conjectured in \eqref{scaling conj}. Nevertheless, robust superdiffusivity in $d=1$ and marginal superdiffusivity in $d=2$ is at least established. Note also, that in $d=1$, in particular cases (of self-repelling lattice walks) the scaling $t^{2/3}X(t)$ has been rigorously established, cf. \cite{toth-17}, \cite{toth-werner-98}, \cite{toth-veto-09}. 

The proof of Theorem \ref{thm:srbp superdiff} follows the \emph{resolvent method} of \cite{landim-olla-yau-98}, \cite{komorowski-olla-02}, \cite{landim-quastel-salmhofer-yau-04}, with new input in the variational computations for the 2-dimensional case. 

Due to the martingale decomposition \eqref{martingale decomposition 2} and stationarity of the process $\eta(s)$, applying a straightforward Schwarz inequality we obtain
\begin{align}
\notag
\int_0^t (t-s) \langle\varphi, P_t \varphi\rangle ds
-
2 \alpha^2 t
\le
\expect{\abs{X(t)}^2}
\le
4 \int_0^t (t-s) \langle\varphi, P_t \varphi\rangle  ds
+
2 \alpha^2 t, 
\end{align}
where $\alpha^2:= t^{-1}\expect{M(t)^2}$ is the variance rate of the first term, $M(t)$,  in the decomposition \eqref{martingale decomposition 2}, which is a square integrable martingale with stationary increments. (The value of $\alpha^2$ is explicitly computable but does not matter.) Hence, 
\begin{align}
\notag
\lambda^{-2} \langle \varphi, R_\lambda \varphi \rangle
-
2\alpha^2\lambda^{-2}
\le
\wh E(\lambda)
\le
4 \lambda^{-2}  \langle \varphi, R_\lambda \varphi \rangle
+
2\alpha^2\lambda^{-2},
\end{align}
where $R_\lambda$ is the resolvent of the semigroup $P_t$. Thus, lower and upper bounds on $\wh E(\lambda)$ reduce to lower and upper bounds on $ \langle \varphi, R_\lambda \varphi \rangle$. The following variational formula, proved in \cite{landim-olla-yau-98},  is valid in the widest generality for any contraction semigroup $P_t=e^{tL}$, with infinitesimal generator $L=-S+A$: 
\begin{align}
\label{variational formula}
\langle \varphi, R_\lambda \varphi \rangle
=
\sup_{\psi\in\cH}
\left\{
2 \langle \varphi,\psi\rangle - \langle\psi, (\lambda I + S) \psi \rangle -\langle A \psi, (\lambda I + S)^{-1} A \psi\rangle
\right\}.
\end{align}
The upper bounds in Theorem \ref{thm:srbp superdiff} are obtained simply by dropping the third (negative!) term on the right hand side of \eqref{variational formula}. This is essentially for free. The lower bounds are obtained by bounding from below the variational expression on the right hand side of \eqref{variational formula}, \emph{in the subspace} $\cH_1$. This leads to a nontrivial variational problem in $u:\R^d\mapsto\R^d$ ($d=1,2$). In $d=2$ the solution is tricky. For details see \cite{toth-valko-12}.

\subsubsection{Diffusive limit in $d\ge3$}
\label{sss: Diffusive limit in d>=3}

\begin{theorem}
\label{thm:srbp clt}
{\rm (Source: \cite{horvath-toth-veto-12b})}
In $d\ge3$, the asymptotic covariance matrix
\begin{align}
\notag
(\sigma^2)_{ij}
:=
\lim_{t\to\infty} t^{-1}\expect{X_i(t)X_j(t)}
\end{align}
exists, it is bounded and non-degenerate.
For any bounded and continuous function $f:\R^d\to\R$, 
\begin{align}
\label{srbp clt in prob}
\lim_{T\to\infty}
\int_\Omega
\abs{
\expectom{f(T^{-1/2}X(T))}
-
\int_{\R^d} \frac{ e^{-\frac {\abs{y}^2} {2} }}{(2\pi)^{\frac{d}{2}}} f(\sigma^{-1}y) dy
}
d\pi(\omega)=0.
\end{align}
\end{theorem}

Theorem \ref{thm:srbp clt} is proved in \cite{horvath-toth-veto-12b}, and weak convergence in the sense of \eqref{srbp clt in prob} of all finite dimensional marginal distributions of the diffusively scaled process $t\mapsto T^{-\frac12} X(Tt)$, as $T\to\infty$, to those of a $d$-dimensional Brownian motion is established.  The proof relies on the efficient martingale approximation \`a la Kipnis-Varadhan of the integral term $I(t)$ on the right hand side of \eqref{martingale decomposition 2}. This is done by verifying the \emph{graded sector condition} of \cite{sethuraman-varadhan-yau-00}. The graded structure of the Hilbert space $\cH$ and of the infinitesimal generator $L$, cf. \eqref{infgen srbp}. Technical details to be found in \cite{horvath-toth-veto-12b}.

\bigskip
\noindent
{\bf Acknowledgement:}
I thank Tomasz Komorowski and Stefano Olla for their help in clarifying some points related to the historical backgrounds.

\hfill\hfill

\hbox{
\hskip55mm
\vbox{\hsize=15cm\noindent
\small
{\sc B\'alint T\'oth}
\\
School of Mathematics, 
University of Bristol
\\
Bristol, BS8 1TW,
United Kingdom
\\
email: {\tt balint.toth@bristol.ac.uk}
}
}

\end{document}